\documentclass{article}

\usepackage{amssymb,dsfont,latexsym}
\usepackage{amsmath}

\newtheorem{thm}{Theorem}

\newtheorem{prop}[thm]{Proposition}
\newtheorem{lem}[thm]{Lemma}
\newtheorem{cor}[thm]{Corollary}

\newcommand\enu[1]{\smallskip\newline\makebox[5mm][l]{\rm(#1)}}

\newcommand\bp{\noindent{\it Proof.}\ }

\newcommand\Ad{{\rm Ad}\,}

\newcommand\rank{{\rm rank}\,}

\author{Erling St{\o}rmer}

\date{07-03-2017}

\title{Mapping cones and separable states}

\begin{document}

\maketitle. 

\begin{abstract}
We study mapping cones and their dual cones of positive maps of the $n\times n$ matrices into itself.  For a natural class of cones there is a close relationship between maps in the cone, super-positive maps, and separable states.  In particular the composition of a map from the cone with a map in the dual cone is super-positive, and so the natural state it defines is separable.
\end{abstract}

\section{Introduction}
One approach to the study of positive maps between operator algebras is via mapping cones.  They were introduced in \cite{S1} and are closed convex cones of positive maps closed under composition with completely positive maps. They are especially useful in finite dimensions as they yield much information on positive maps with given positivity properties, see e.g. \cite{S}.  In the present paper we shall study mapping cones which are closed under composition with all positive maps.  Then we obtain characterizations of super-positive maps, also called entanglement breaking maps, and their relationship to separable operators and states.

Three classes of maps will be central in the present paper.  Let $M_n$ denote the complex $n\times n$. matrices, and let  $\phi \colon M_n  \to M_n$ be a linear map. Then $\phi$ is \textit{positive}, written $\phi\geq 0$, if $\phi(a) \geq 0$ whenever $ a \in M_n^+$, the positive matrices in $M_n$.    $\phi$ is \textit{completely positive} if
 $ \iota \otimes \phi \geq 0$, where $\iota$ is the identity map on $M_n$.  Then $\phi$ is of the form $\phi = \sum_i AdV_i$, finite sum, where $V_i \in M_n$, and $AdV(a) = V^* aV$, see \cite{S} Thm. 4.1.8.  $\phi $ is \textit{super-positive} if it is of the form $\sum_i a_i \omega_i$, where $a_i \in M_n^+ $, and $\omega_i$ is a state on $M_n$. 
 The study of positive maps $\phi$ is equivalent to the study of its \textit{Choi matrix} denoted by  $C_\phi$ and defined as follows.
 $$ 
 C_\phi = \sum  e_{ij} \otimes \phi(e_{ij}) \in M_n \otimes M_n,
 $$
where $e_{ij}, 1\leq i,j \leq n$, is a complete set of matrix units for $M_n$. If $Tr$ is the usual trace on $M_n \otimes M_n$ then we have, see \cite {S}, Ch. 4.
$\phi \geq 0$  iff  $ Tr(C_{\phi} a\otimes b) \geq 0$    for all $a, b \in M_n^+$ , and this holds   iff $Tr(C_{\phi} C_{\psi} )\geq 0$  for all super-positive $\psi$.

$\phi$ is completely positive iff $C_{\phi} \geq 0$.  $\phi$ is super-positive iff $Tr(C_{\phi}  C_{\psi}) \geq 0$ for all positive maps $\psi$ iff $C_{\phi}$ is separable, i.e.
$$
C_{\phi} = \sum a_i \otimes b_i,  a_i, b_i \in M_n^+.
$$
We shall use the notation

$P$ is the cone of all positive maps of $M_n$ into itself. 

$CP$ is the cone of all completely positive maps.

$SP$ is the cone of all super-positive maps.

\textit{Acknowledgement.}The author is indebted to Geir Dahl for helpful comments on  dual cones.
 
\section{Mapping cones and their dual cones.}

As above $P$ denotes the set of positive maps of $M_n$ into itself.  A \textit{mapping cone} is a closed convex cone $J \subset P$ such that $\phi \in J$ implies
$\alpha \circ \phi \circ \beta \in J$ whenever $\alpha, \beta \in CP$.  We always assume $J$ is \textit{symmetric}, i.e. $\phi \in J$ implies $\phi^* \in J$, and $\phi^t \in J$, where $\phi^*$ and $\phi^t$ are defined by
$
Tr(\phi(a)b) = Tr(a\phi^*(b)),   \phi^t(a)= \phi(a^t)^t, \forall  a,b \in M_n,
$
where t denotes the transpose.

$J$ is an \textit{invariant mapping cone} if $\phi \in J$ implies $\alpha \circ \phi \circ \beta \in J$ for all $\alpha, \beta \in P$

The \textit{dual cone} of $J$ is the cone
$$
J^o = \{ \phi \in P: Tr(C_{\phi} C_{\psi}) \geq 0, \forall \psi \in J\}.
$$
Then $J^{oo} = J$, see \cite{S}, Section 6.2, and $J^o$ is a mapping cone which is invariant if $J$ is. The definitions of $P,CP,SP$ in the introduction can be restated as
\begin{eqnarray*}
\phi \in P \Leftrightarrow \phi \in SP^o.\\
\phi \in CP \Leftrightarrow \phi \in CP^o.\\
\phi \in  SP \Leftrightarrow  \phi \in P^o.
\end{eqnarray*}
The next result with a different proof is the same as \cite{S}, Lem.5.1.5.

\begin{lem}\label{lem1}
$SP$ is an invariant mapping cone, and is also a minimal mapping cone.
\end{lem}
\bp  $SP$ is clearly an invariant mapping con since if $\omega$ is a state and $\phi \in P$ then $\omega \circ \phi$ is a scalar multiple of a state, so belongs to $SP$, and 
$\phi  \circ \omega(a) = \omega(a)\phi(1) $ is also in $SP$.  If $\psi \in SP$ and  $\phi \in P$, then $Tr(C_{\phi} C_{\psi}) \geq 0.$ In particular this holds for $\phi \in J^o$, hence 
$\psi \in J^{oo} = J$.  Thus $SP \subset J$, completing the proof.
\medskip

Note that $SP$ is not a hereditary cone, because the trace $Tr_n$ on $M_n$ belongs to $SP$, and by \cite{S}, Thm.  7.5.4, if $\phi \in P$ then $\psi = \phi(1)Tr_n + \phi \in SP$, and so $\phi = \psi - \phi(1)Tr_n$ is the difference of two maps in $SP$, hence $\phi \leq \psi$ even though $\phi$ need not belong to $SP$.

The following theorem is central for our development of the theory to follow.  It is a slight  variation of Theorems 6.6 and 7.1.1 in \cite {S}.  Recall that if $\phi \in P$ then
$\widetilde{\phi}$ is the linear functional on $M_n \otimes M_n$ defined by
$$
\widetilde{\phi}(a \otimes b) = Tr(\phi(a) b^t).
$$
By \cite{S}, Lem. 4.2.3, $C_{\phi}^t$ is the density matrix for $\widetilde\phi$.  Thus $\widetilde \phi$ is positive iff $\phi$ is completely positive.  By \cite{S} Prop. 5.1.4 
$\widetilde \phi$ is separable iff $C_{\phi}$ is separable.

\begin{thm}\label{thm2}
 Let $J$ be aa mapping cone and $\phi \in P$.  Then the following conditions are equivalent.
 \enu{i}  $\phi \in J^o.$
 \enu{ii}  $\phi \circ \alpha \in CP \, \forall \alpha \in J$.
 \enu{iii} $ \alpha \circ \phi \in CP\,    \forall \alpha \in J.$
 \enu{iv} $\widetilde\phi \circ (\iota \otimes \alpha)\geq 0\,   \forall \alpha \in J$.
 \enu{v}$(\iota \otimes \alpha)(C_{\phi}) \geq 0 \,  \forall  \alpha \in J.$
 \end{thm}
\bp The only part of the theorem which is not contained in \cite{S},Thm. 6.1.6 is the equivalence of (ii) with the others. But we have by symmetry of $CP $ that
$$
\phi \circ \alpha \in CP\Leftrightarrow (\phi \circ \alpha)^* = \alpha^* \circ \phi^* \in  CP.
$$
Since $J$ is symmetric, as we have assumed for all our mapping cones, $\phi\circ \alpha \in CP $ iff  $ \alpha \circ \phi^* \in CP$ for all $\alpha \in J$, hence by the equivalence (i) $\Leftrightarrow $(iii), iff  $\phi^* \in J^o.$  But $J^o$  is symmetric, so $\phi^* \in J^o$ iff $\phi\in J^o$, hence we have (i) $\Leftrightarrow$ (ii), completing the proof.

If $J$ is an invariant mapping cone and $\phi \in P$, conditions (ii) and (iii) can be sharpened.

\begin{cor}\label{cor3}
Let $J$ be an invariant mapping cone and $\phi \in P$.  Then the  conditions in Theorem 2 are equivalent to
\enu{vi} $\phi \circ \alpha \in SP\, \forall \alpha \in J.$
\enu{vii} $\alpha \circ \phi \in SP\,  \forall \alpha \in J.$
\end{cor} 
\bp
Let $\alpha, \psi \in P$.  Then we have for $\alpha \in J$
$$
Tr(C_{\alpha\circ\phi} C_{\psi}) = Tr( \iota \otimes \alpha( C_{\phi})C_{\psi}) = Tr(C_{\phi} C_{\alpha^* \circ \psi}).
$$
Hence 
$$
Tr(C_{\alpha \circ \phi} C_{\psi}) \geq 0 \, \forall \psi\, \Leftrightarrow \phi \in J^o,
$$
as the maps $\alpha^* \circ \psi$ include all maps in $J$. Thus $\phi \in J^o$ iff $\alpha \circ \phi \in SP $ for all $\alpha \in J$. As in the proof of the equivalence (ii) $\Leftrightarrow$ (iii) in Theorem 2 we show (vi) $\Leftrightarrow $ (vii).  The proof is complete.

\medskip
If $V$ is a finite dimensional Hilbert space, and $J$ and $K$ are closed convex cones in $V$, then the dual cone $(J \cap K)^o$ of their intersection is the closure of the sum $J^ o + K^o$ of their dual cones, see \cite{B}. From this we get the following corollary, which is included for completeness.

\begin{cor}\label{cor4}
Let $J$ and $K$ be mapping cones in $P$.  Then $(J \cap K)^o  = J^o + K^o$.
\end{cor}
\bp
By the above and finite dimensionality it suffices to show that $J^o + K^o$ is closed in the norm topology.  So let $(\phi_i = \alpha_i + \beta_i)$ be a sequence converging to
$\phi$, where $\alpha_i \in J^o, \beta_i \in K^o$. We can assume $\parallel \phi_i \parallel  \leq \parallel \phi \parallel$ for all $i$ and the same for $\alpha_i$ and $\beta_i$.  Since the dimension is finite., by compactness of the unit ball in $J^o$ there is a subsequence of $(\alpha_i)$ which converges to a map $\alpha \in J^o$ and a subsequence $(\beta_j)$ of $(\beta_i)$ which converges to $\beta \in J^o$.  Thus 
$$
\phi_j = \alpha_j + \beta_j  \longrightarrow \phi = \alpha + \beta \in J^o + K^o,
$$
proving the corollary.
\medskip 

\textit{Remark}. In the sequel we shall see that maps of the form $Ad e$ defined by $Ad e(a) = eae$ are of central importance. If $a \in M_n$ with range  $range(a) =e$, then for $J$ a mapping cone $Ad e \in J$ iff $Ad a \in J$.  Indeed, we have $a = ae$, so $ Ad a = Ad a \circ Ad e \in J$ if $Ad e \in J$, and if $Ad. a \in J,$ since $range(aa^*) = e$ we can assume $a \geq 0$.  Let $a^{-1}$ be the inverse of $a$ in $eM_n e$, so $aa^{-1} = e$.  $Ad e = Ad a \circ Ad a^{-1} \in J$, proving the converse. 

\begin{thm}\label{thm5}
Let $J$ be an invariant mapping cone.  Let
$$
%:
k = max \{ \rank( e):  Ad e \in J \}.
$$
Assume $k \geq 2. $. Then we have
\enu{i} If  $\phi \in J^o$ and $rank(range (\phi)) \leq k$, then $\phi \in SP$.
\enu{ii} If $e$ is a projection in $M_n$ with  $rank(e) \leq k$, then $Ad e \in J$.

\end{thm}
\bp
Let $\phi \in J^o$ be as in (i).  By assumption there exists a projection $e$ of rank k such that $Ad e \in J$.  Let $f $ be the range projection of $\phi$, so by assumption 
 $rank(f) \leq rank(e)$.  Thus there exists a partial isometry $v \in M_n$ such that $f = v^* e v$.  Hence
 $$
 Ad f = Ad v \circ Ad e \circ Ad v^* \in J.
 $$
 Thus $\phi = Ad f \circ \phi \in J$, hence $\phi$ is the composition of a map in $J$ with a map in $J^o$, so $\phi \in SP$ by Corollary 3, proving (i).
 To show (ii) let $f$ be a projection of rank k such that $f\geq e$. By the above argument $Ad f \in J$ hence $Ad e = Ad f \circ Ad e \in J$.

\begin{cor}\label{cor6}  
Let n=2 and $J$ be a an invariant mapping cone in $P$.  Then either $J = P$ or $J = SP$.
\end{cor}
\bp
Let $k$ be as in Theorem 5.  Since every map in $P$ is decomposable, see \cite{S},Thm. 6.3.1, each map is a linear combination of maps of the form $Ad v$ or $Ad v \circ t,$ where    $v \in M_2$, and  $t$ denotes the transpose.  If $k=1$ it thus follows  that all maps in $J$ are super-positive, hence in $SP$.  If $k=2$ then the identity map belongs to $J$, so $J = P$, completing the proof of the corollary.

\medskip
We denote by $D_2$ the invariant mapping cone generated by $Ade$ for the rank 2 projections in $M_n$, or by the remark before Theorem 5, matrices $a$ of rank 2. By the argument in the proof of Theorem 5 $D_2$ is the invariant mapping cone generated by a single map $Ad e$ with e of rank 2.  By Theorem 5, if $n>2$, and Corollary 6 if $n=2$, we have

\begin{cor}\label{cor7}
If $e$ is a projection with $Ad e \in D_2^o$, and $rank(e). \leq 2$ then $rank(e) = 1$.
\end{cor}

\begin{thm}\label{thm8}
Let $J$ be an invariant mapping cone not contained in $D_2^o$. Then $J $ contains $D_2$.
\end{thm}
\bp
Let $\phi \in J$,  but not in $D_2^o$.  Then in particular $\phi$ does not belong to $SP$.  We assert that there exist $\alpha, \beta, \gamma \in P$ and projections $e, f$ of rank 2 such that 
$$
Ad e \circ \alpha \circ \phi \circ \beta \circ Ad f \circ \gamma \notin SP.
$$
If the assertion is false, i.e. the above map belongs to $SP$ for all $\alpha, \beta, \gamma, e,f$, then, since the maps $\beta \circ Ad f \circ \gamma$ generate $D_2$, it follows from Corollary 3 that
$$
Ad e \circ \alpha \circ \phi \in D_2^o \, \, \forall e, \alpha.
$$
Again, since $Ad e \in D_2$, by Corollary 3,
$$
Ad e \circ \alpha \circ \phi = Ad e \circ (Ad e \circ \alpha \circ \phi)  \in SP.
$$
It follows that from the last equation
$$
(\eta \circ  Ad e \circ \alpha) \circ \phi \in SP\, \,  \forall \eta, \alpha \in P.
$$
Since the maps $\eta \circ Ad e \circ \alpha$ generate $D_2$, it follows again from Corollary 3 that $\phi \in D_2^o$, contrary to assertion at the beginning of the proof. Thus the assertion follows. We can therefore find projections $e$ and $f$ of rank 2 and $\alpha, \beta \in P$ such that$$
\Ad e \circ \alpha \circ \phi \circ \beta \circ Ad f  \notin SP.
$$
As in the proof of Theorem 5 (i) there exists a partial isometry $v$ such that $Ad f = Adv \circ Ad e \circ Adv^*$, so we can replace $f$ by $e$, and thus find $e$ of rank 2 such that
$$
Ad e \circ (\alpha \circ \phi \circ \beta)\circ Ad e \notin SP.
$$
Let
$$
J_2 = \{ Ad e \circ \psi \circ Ad e: \psi \in J \}.
$$
Since $\alpha \circ \phi \circ \beta \in J$ it follows by the above that $J_2 \nsubseteq SP$. $J_2$ is an invariant mapping cone inside the set 
$$
\{ Ad e \circ \psi \circ Ad e: \psi \in P\},
$$
which is linearly isomorphic to the positive maps of $M_2$ into itself. By Corollary 6, since $J_2 \nsubseteq SP$ 
$$
J_2 = \{ Ad e \circ \psi \circ  Ad e: \psi \in P \}.
$$
But then $Ad e \in J_2$  By assumption $\phi \in J$, and $J_2 \subset J$, so $Ad e \in J$.  Since $Ad e$ generates $D_2$ as an invariant mapping cone, we thus conclude that $D_2 \subset J$. The proof is complete
\medskip 

A map $\phi \in P$ is called \textit{extremal} if whenever $\psi \in P$ and $\phi \geq \psi$, then $\psi = \lambda \phi$ for some $0\leq \lambda \leq 1$.  This can be formulated as, if $\phi = \sum_{i=1} ^k \psi_i$   is a sum of positive maps $\psi_i$, then $k =1$, so $\phi = \psi_1$.  It is shown in \cite{S} Prop. 3.1.3. that each map of the form $\phi = Ad v$ is extremal.

\begin{prop}\label{prop 9}.
Let $\phi\in P$ be an extremal map, and let $J$ be the invariant mapping cone generated by $\phi$.  Then we have:
\enu{i} If $\phi \notin D_2\bigvee D_2^o$ then $J \supsetneqq D_2$, and $J \nsubseteq D_2^o.$
\enu{ii} If $\phi$ is unital then $J \nsubseteq P$ unless $\phi$ is either an automorphism or anti-automorphsm.
\end{prop}
\bp
Ad (i).  Clearly $J \nsubseteq D_2^o$, hence by Theorem 8 $ J \supseteq D_2$. Since $\phi \notin D_2, J\nsupseteq D_2.$
Ad(ii). If $J = P$ then the identity map $\iota\in J.$   Since $\iota$ is an extremal map it must be of the form $\iota = \alpha \circ \phi \circ \beta$ with $\alpha,\beta \in P$.  Furthermore, $\alpha$ and $\beta$ must be extremal, since otherwise $\iota$ would be a sum of several positive maps.  We thus have $\beta$ is invertible with inverse $\alpha \circ \phi$.  Thus $\beta$ is an order-isomorphism of $M_n$ onto itself.  Similarly $\alpha$ is an order-isomorphism, so  that $\phi$ is the same.  Since $\phi$ is unital, $\phi$ is either an automorphism or an anti-automorphism, see e.g. \cite{S} Thm.2.1.3. The proof is complete.

\section{PPT-states}
\textit{PPT-states} were introduced by Peres \cite{P} in 1996 and are states $\rho$ on $M_n \otimes M_n$ such that $\rho\circ (\iota \otimes t)$ is a state, where $t$ as before is the transpose on $M_n$.  For a while it was believed that they were separable, but a counter example was exhibited by P. Horodecki \cite {PH} in 1997.

We call a map $\phi \in P$ for a \textit{PPT-map} if the corresponding linear functional $\widetilde\phi$ is a PPT-state. The following consequences can be read out of \cite{S}, Section 7.2. 
\begin{prop}\label{prop10}
Let $\phi \in P$.  Then the following conditions are equivalent.
\enu{i}  $\phi \in CP \cap coCP$.
\enu{ii} $\phi \circ \alpha \in CP \, \, \forall \alpha \in CP \bigvee coCP.$
\enu{iii} $\phi$ is a PPT-map.
\end {prop}

Recall that map $\phi \in CP \bigvee coCP$  is decomposable, so it is of the form $\phi =  \alpha + \beta$ with $\alpha \in CP,  \,  \beta \in coCP$, see Corollary 4.  The next result is a dual version of Proposition 10.  In (ii) we used that $CP \cap coCP$ is a mapping cone with dual cone $CP \bigvee coCP$.

\begin{prop}\label{prop11}.
Let $\phi \in P$.  Then the following conditions are equivalent.
\enu{i} $ \phi \in CP \bigvee coCP$.
\enu{ii} $\phi \circ \alpha \ \in CP \cap coCP, \, \,\forall \alpha \in CP\bigvee coCP.$
\enu{iii} $ \widetilde\rho(C_{\phi}) \geq 0,  \, \forall  PPT$-maps $\rho$
\end{prop}
If n=3, then $D_2$ consists of decomposable maps.  For completeness we include a proof of this fact.

\begin{lem}\label{lem12}
Let n=3. Then $D_2$ consists of decomposable maps, hence $D_2 \subset CP \bigvee coCP$.
\end{lem}
\bp
Let $\phi\in D_2$.  Then $\phi$ is a sum of maps $\alpha \circ Ade \circ \beta$ with $e$ a projection of rank $\leq 2$, and $\alpha, \beta \in P$.  Each map $\alpha \circ Ade$ can be identified with a map $M_2 \longrightarrow M_3$, and the map $Ade \circ \beta$ with a map $M_3\longrightarrow M_2$. By a result of Woronowicz \cite {W}, both maps are decomposable, hence so is $\alpha \circ Ade \circ\beta$, and therefore $\phi$ is decomposable, , so belongs to $CP\bigvee coCP$, proving the lemma.
\medskip 

\begin{prop}\label{prop13}.
Let n=3 and $\phi \in P$ be a PPT-map. Then we have:
\enu{i} $\phi \in D_2^o$.
\enu{ii} $\phi \circ \alpha \in SP$ and $\alpha \circ \phi \in SP\, \forall \alpha \in D_2$.
\enu{iii} If $range(\phi)$ or $supp(\phi)$ is a projection of rank 2, then $\phi \in SP$
\end{prop}
\bp
By Proposition 10 and Lemma 12  $\phi \in CP \cap coCP = (CP \bigvee coCP)^o \subset D_2^o$, proving (i). Thus by  Corollary 3, $\phi \circ \alpha \in SP$, and
 $\alpha \circ \phi \in SP$ for all $\alpha \in D_2$, proving (ii).
 To show (iii) note that under the assumptions on range and support of $\phi,$   $Ad (range(\phi))$ or  $Ad (supp(\phi))$ belongs to $D_2$.  Hence, either  
 $\phi = Ad(range(\phi))\circ\phi \in SP$ by (ii), or similarly for $\phi = \phi \circ Ad(supp(\phi)) \in SP$, so $ \phi \in SP$ completing the proof.
 \medskip 
 
 In \cite{MH} it was shown that all PPT-states on $M_2\otimes M_2$ and $M_3\otimes M_2$ are separable.  These results follow easily from the above results.  Indeed if n=2 then by Proposition 10 and Corollary 6 we get
 $$
 \phi \in CP \cap coCP = (CP \bigvee coCP)^ o \subset D_2^o = SP.
 $$
 When n=3, and $e = range(\phi)$ is of rank 2 then, since $eM_3 e \simeq M_2$, we can consider  $\phi$ as a PPT-map from $M_3$ to $M_2$. By Proposition 13(iii) $\phi \in SP$. 
  To translate  this to states recall that $C_{\phi}^t$ is the density matrix for $\widetilde\phi$, and $\widetilde \phi$ is separable iff $\phi \in SP$.

Department of Mathematics,  University of Oslo, 0316 Oslo, Norway.

e-mail  erlings@math.uio.no

\end{document}